\documentclass[12pt,a4paper]{article}
\usepackage{mathrsfs}

\usepackage{amssymb}

\usepackage{amsmath}

\newtheorem{theorem}{Theorem}[section]
\newtheorem{definition}[theorem]{Definition}
\newtheorem{lemma}[theorem]{Lemma}
\newtheorem{corollary}[theorem]{Corollary}

\date{}

\begin{document}

\title{Gr\"obner--Shirshov bases for Vinberg--Koszul--Gerstenhaber
right-symmetric algebras \footnote{Supported by the NNSF of China
(No.10771077) and the NSF of Guangdong Province (No.06025062).}}

\author{
L. A. Bokut\footnote {Supported by the grant LSS--344.2008.1 and SB
RAS
Integration grant No. 2009.97 (Russia).} \\
{\small \ School of Mathematical Sciences, South China Normal
University}\\
{\small Guangzhou 510631, P. R. China}\\
{\small Sobolev Institute of Mathematics, Russian Academy of
Sciences}\\
{\small Siberian Branch, Novosibirsk 630090, Russia}\\
{\small Email: bokut@math.nsc.ru}\\
\\
 Yuqun
Chen\footnote {Corresponding author.} \  and Yu Li\\
{\small \ School of Mathematical Sciences, South China Normal
University}\\
{\small Guangzhou 510631, P. R. China}\\
{\small Email: yqchen@scnu.edu.cn}\\
{\small liyu820615@126.com}}

\maketitle \noindent\textbf{Abstract:} In this paper, we establish
the Composition-Diamond lemma for right-symmetric algebras. As an
application, we give a Gr\"{o}bner-Shirshov basis for universal
enveloping right--symmetric algebra of a Lie algebra.



\section{Introduction}

 Gr\"{o}bner  and  Gr\"{o}bner--Shirshov bases theories were
 invented independently by A.I. Shirshov
\cite{Sh62c} for Lie algebras and H. Hironaka \cite{Hi64} and B.
Buchberger \cite{Bu65,Bu70} for associative-commutative algebras.

Shirshov's paper \cite{Sh62c} based on his papers \cite{Sh62b}
(Gr\"{o}bner--Shirshov bases theory for (anti) commutative algebras,
the reduction algorithm for (anti-) commutative algebras) and
\cite{Sh58} (Lyndon--Shirshov words which were defined some earlier
in \cite{Lyndon54}, but incidentally that was unknown for 25 years
in Russia and these words were called Shirshov's regular words, see,
for example, \cite{Bah87,Bah03,BMPZ92,Bo72,Ufn95,Mik96}, see also
\cite{Cohn65}), Lyndon--Shirshov basis of a free Lie algebra (see
also \cite{CFL58}). The latter Shirshov's papers \cite{Sh62b,Sh58}
were based on his Thesis \cite{Sh53a}, A.G. Kurosh (adv), published
in three papers \cite{Sh53b} (on free Lie algebras: $K_d$-Lemma
(Lazard--Shirshov elimination process), the subalgebra theorem
(Shirshov--Witt theorem)), \cite{Sh54} (on free (anti-) commutative
algebras: linear bases, the subalgebra theorems), and \cite{Sh62a}
(on free Lie algebras: a series of bases with well-ordering of basic
Lie words such that $[w]=[[u][v]]>[v]$, see also \cite{Vie78}; the
series is called now Hall sets \cite{Reu93} or Hall--Shirshov
bases). Shirshov's Thesis, in turn, was in line with a Kurosh's
paper \cite{Ku60} (on free non-associative algebras: the subalgebra
theorem). Also Shirshov's paper \cite{Sh62b} was in a sense of a
continuation of a paper by A.I. Zhukov \cite{Zhu50}, a student of
Kurosh (on free non-commutative algebras: decidability of the word
problem for non-associative algebras). The difference from the
Zhukov's approach was that Zhukov did not use any linear ordering of
non-associative words, but just the partial deg-ordering to compare
two words by the degree (length).

\ \

It would be not a big exaggeration to say that Shirshov's paper
\cite{Sh62c} was between line of the Kurosh's program of study free
algebras of different classes of non-associative algebras.

\ \

Shirshov's paper \cite{Sh62c} contained implicitly the
Gr\"{o}bner--Shirshov bases theory for associative algebras too
because he constantly used that any Lie polynomial is at the same
time a non-commutative polynomial. For example, the maximal term of
a Lie polynomial is defined as its maximal word (in the deg-lex
ordering) as a non-commutative polynomial, the definition of a Lie
composition (Lie S-polynomial) of two Lie polynomials begins with
the definition of their composition as non-commutative polynomials
and follows by putting some special Lie brackets on it, and so on.
The main Composition (-Diamond) lemma for associative polynomials is
actually proved in the paper and we need only to ``forget" about Lie
brackets in the proof of this lemma for Lie polynomials, i.e. to
change Lie polynomials to non-commutative (associative) polynomials
(\cite{Sh62c}, Lemma 3). Explicitly Composition (-Diamond) lemma was
formulated much later in papers L.A. Bokut \cite{Bo76} and G.
Bergman \cite{Be78}.

Last years there were quite a few results on Gr\"{o}bner--Shirshov
bases for associative algebras, Lie (super) algebras and irreducible
modules for them, Kac--Moody algebras, Coxeter groups, braid groups,
quantum groups, conformal algebras, free inverse semigroups,
Kurosh's $\Omega$-algebras, Loday's dialgebras and Leibniz algebras,
Rota-Baxter algebras, and so on, see, for example, books
\cite{BoK94, MZ95}, papers \cite{MZ95, La97, BoKLM99, La00a, La00b,
La05,  KL00a,KL00b,
BokutShiao01,KLLO02,BChainikovShum07,Bo07,KLLP07, BFKS08, BFK04,
Bo09}, and surveys \cite{BFKK00, BK00, BK05a,
BChibrikov05,BokutChen08, {BCZ08}}. Actually, conformal algebras,
dialgebras, Rota-Baxter algebras are examples of $\Omega$-algebras.
For non-associative $\Omega$-algebras, Composition-Diamond lemma was
proved in \cite{DH08}. The case of associative $\Omega$-algebras
(associative algebras with any set $\Omega$ of multi-linear
operations) was treated in \cite{BCQ08} with an application to free
($\lambda$-) Rota-Baxter algebras (the latter are associative
algebras with linear operation $P(x)$ and the identity
$P(x)P(y)=P(P(x)y) + P(xP(y)) +\lambda P(xy)$, where $\lambda$ is a
fix element of a ground field, see, for example, \cite{EGuo08}).
Composition--Diamond lemma for dialgebras \cite{BCL08} has an
application to the PBW theorem for universal enveloping dialgebras
of  Leibniz algebras (see \cite{Aymon03}).

In this paper, we are doing Gr\"obner-Shirshov bases theory for
right-symmetric algebras (RS-algebras) (they are defined by the
following identity $(x,y,z)=(x,z,y)$ for the associator
$(x,y,z)=(xy)z-x(yz) $. Under this name, these algebras appeared in
the paper E.B. Vinberg (\cite{Vin60}, 1960) (actually, he invented
left-symmetric (LS-) algebras, $(x,y,z)=(y,x,z)$). Independently
they were introduced by J.-L. Koszul (\cite{Ko61}, 1961) and M.
Gerstenhaber (\cite{Ge63}, 1963) (under the name pre-Lie algebras).
As was pointed out by D. Burde \cite{Bur06}, the same algebras first
appeared in a paper by A. Cayley in 1896 (see \cite{Cay}). Survey
\cite{Bur06} contains detailed discussion of the origin, theory and
applications of LSA in geometry and physics together with an
extensive bibliography.

D. Segal \cite{Se94} found a linear basis of a free LS-algebra and
applied it for PBW type theorem for universal LS enveloping algebra
of a Lie algebra. Vasil'eva and Mikhalev \cite{VM96} found another
proof of the former Segal's result (and a more general result for LS
superalgebras) using Composition-Diamond lemma from Zhukov's paper
\cite{Zhu50}. Kazybaev, Makar-Limanov and Umirbaev \cite{KMLU08}
found some new properties of the Segal's basic words and proved the
Freiheitssatz and decidability of the word problem for one-relator
RS-algebras. Actually, it is a generalization of well known
Shirshov's results for one-relator Lie algebras \cite{Sh62b}.

\section{Composition-Diamond lemma for right-symmetric algebras}

Let $X=\{x_i|i\in I \}$ be a set, $X^*$ the set of all associative
words $u$ in $X$, $X^{\ast \ast }$ the set of all non-associative
words $(u)$ in $X$, and $|(u)|$ the length of the word $(u)$.

Let $I$ be a well ordered set. We order $X^{**}$ by the induction on
the lengths of the words $(u)$ and $(v)$ in $X^{**}$:
\begin{enumerate}
\item[(i)] \  If $|((u)(v))|=2$, then $(u)=x_i > (v)=x_j$ if and
only if $i>j$.
\item[(ii)] \ If $|((u)(v))|>2$, then $(u)>(v)$ if and
only if one of the following cases holds:
\begin{enumerate}
\item[(a)] $|(u)|>|(v)|$.
\item[(b)] \ If $|(u)|=|(v)|, \ (u)=((u_1)(u_2))$ and
$(v)=((v_1)(v_2))$, then $(u_1)>(v_1)$ or ($(u_1)=(v_1)$  and
$(u_2)>(v_2)$).
\end{enumerate}
\end{enumerate}

It is clear that the order $<$ on $X^{**}$ is well ordered. This
order is called deg-lex (degree-lexicographical) order and we use
this order throughout this paper.

We now cite the definition of good words (cf. \cite{Se94}) by
induction on  length:
\begin{enumerate}
\item[1)]$x_i$ is a good word for any $x_i\in X$.

Suppose that we define good words of length $<n$.

\item[2)] non-associative word $((v)(w))$ is called a good word if
\begin{enumerate}

\item[(a)] both $(v)$ and $(w)$ are good words,

\item[(b)] if $(v)=((v_1)(v_2))$, then $(v_2)\leq(w)$.
\end{enumerate}
\end{enumerate}

We denote $(u)$ by $[u]$, if $(u)$ is a good word.

Let $W$ be the set of all good words in the alphabet $X$ and
$RS\langle X\rangle$ the free right-symmetric algebra over a field
$k$ generated by $X$. Then $W$ forms a linear basis of the free
right-symmetric algebra $RS\langle X\rangle$, see \cite{Se94}.

Every nonzero element $f$ in $RS\langle X\rangle$ can be uniquely
represented as
$$
f=\lambda_1[w_1]+\lambda_2[w_2]+\dots+\lambda_n[w_n]
$$
where $[w_i]\in W$, $0\neq\lambda_i\in k$ for all $i$ and
$[w_1]>[w_2]>\dots>[w_n]$. Denote by $\bar{f}$ the leading word
$[w_1]$ of $f$. $f$ is called monic if the coefficient of $\bar{f}$
is 1.

For any $(w),(w_1)\in X^{**}$, denote by
$$(w)R_{(w_1)}=((w)(w_1)).
$$
The following results are actually proved in \cite{KMLU08}.

\begin{lemma}\label{l1}(\cite{KMLU08})
In $X^{**}$, every good word $[w]\in W$  can be uniquely represented
as
$$
[w]=x_iR_{[w_1]}R_{[w_2]}\cdots R_{[w_n]}
$$
where $[w_j]\in W$ for all $j$ and $[w_1]\leq{[w_2]}\leq\cdots \leq
{[w_n]}$.

\end{lemma}

\begin{lemma}\label{l3}(\cite{KMLU08})
Let $[u]$ and $[v]$ be arbitrary good words and assume that
$[u]=\newline x_iR_{[u_1]}R_{[u_2]}\cdots R_{[u_n]}$. Then, in
$RS\langle X\rangle$,
$$
\overline{[u][v]}=x_iR_{[u_1]}\cdots
R_{[u_s]}R_{[v]}R_{[u_{s+1}]}\cdots R_{[u_n]}
$$
where $[u_1]\leq \cdots \leq [u_s]\leq [v]<[u_{s+1}]\leq\cdots \leq
[u_n]$ and $s\leqslant n$.
\end{lemma}

By Lemma \ref{l3}, we have

\begin{corollary}\label{ll5}
Let $[u], \ [v]\in W$ and $[u]=x_iR_{[u_1]}\cdots
R_{[u_{m-1}]}R_{[u_m]}=[u'][u_m]$, where $[u']=x_iR_{[u_1]}\cdots
R_{[u_{m-1}]}$. If $[u_m]>[v]$, then, in $RS\langle X\rangle$,
$\overline{[u][v]}=\overline{([u'][v])}[u_m]$ and
$\overline{[u']([u_m][v])},\newline \overline{[u']([v][u_m])}<
\overline{[u][v]}$.
\end{corollary}

\begin{lemma}\label{l4}(\cite{KMLU08})
Let $[u],[v]$, and $[w]$ be arbitrary good words. If $[u]<[v]$ then
$\overline{[w][u]}< \overline{[w][v]}$ and
$\overline{[u][w]}<\overline{[v][w]}$. It follows that if $f,g\in
RS\langle X\rangle$, then $\overline{fg}=\overline{\bar{f}\bar{g}}$.
\end{lemma}

\begin{definition}
Let $S\subset A$ be a set of monic polynomials, $s\in S$ and $(u)\in
X^{**}$.
We view $S$ as a new set of letters with $S\cap X=\emptyset$. We define $S$-word $%
(u)_s$ by induction:
\begin{enumerate}
\item[(i)] $(u)_s=s$ is an $S$-word of $S$-length
1.
\item[(ii)] If $(u)_s$ is an $S$-word with $S$-length k and $[v]$
is a good word with length $l$, then
$$
(u)_s[v]\ and\ [v](u)_s
$$
are $S$-words with $S$-length $k+l$.
\end{enumerate}
\end{definition}

\begin{definition}
An $S$-word $(u)_s$ is called a normal $S$-word if $(u)_{\bar s
}=(a\bar s b)$ is a good word. We denote $(u)_s$ by $[u]_s$ if
$(u)_s$ is a normal $S$-word. From Lemma \ref{l4}, it follows that
$\overline{[u]_s}=[u]_{\bar s }$.
\end{definition}

\begin{definition}
Let $f,g\in RS\langle X\rangle$ be monic polynomials, $[w]\in W$,
and $a,b\in X^{*}$. Then there are two kinds of compositions.
\begin{enumerate}
\item[(i)] If $\bar{f}=[a\bar{g}b]$, then
$$
(f,g)_{\bar{f}}=f-[agb]
$$
is called composition of inclusion.

\item[(ii)] If $(\bar{f}[w])$ is not good, then
$$
f\cdot [w]
$$
is called composition of right multiplication.
\end{enumerate}
\end{definition}

Let $S\subset RS\langle X\rangle$ be a given  nonempty subset. The
composition of inclusion $(f,g)_{\bar{f}}$ is said to be trivial
modulo $S$ if
\begin{equation*}
(f,g)_{\bar{f}}=\sum\limits_i\alpha_i[a_is_ib_i]
\end{equation*}
where each $\alpha_i\in k,\ a_i,b_i\in X^*,\ s_i\in S,\ [a_is_ib_i]$
is normal S-word and $[a_i\bar{s_i}b_i]<{\bar{f}}$. If this is the
case, then we write
$$
(f,g)_{\bar{f}}\equiv 0\ \  \ mod (S,{\bar{f}}).
$$

In general, for any normal word $[w]$  and $p,q\in RS\langle
X\rangle$, we write
$$
p\equiv q\quad mod(S,[w])
$$
which  means that $p-q=\sum\alpha_i [a_i s_i b_i] $, where each
$\alpha_i\in k, \ a_i,b_i\in X^{*}, \ s_i\in S$ and $[a_i \bar {s_i}
b_i]<[w]$.

The composition of right multiplication $f\cdot [w]$ is said to be
trivial modulo $S$ if
$$
f\cdot [w]=\sum\limits_i\alpha_i[a_is_ib_i]
$$
where each $\alpha_i\in k,\ a_i,b_i\in X^*,\ s_i\in S,\ [a_is_ib_i]$
is normal S-word and $[a_i\bar{s_i}b_i]\leq \overline{f\cdot [w]}$.
If this is the case, then we write
$$
f\cdot [w]\equiv 0\ \  \ mod (S).
$$

\begin{definition}
Let $S\subset RS\langle X\rangle$ be a nonempty set of monic
polynomials and the order $<$ on $X^{**}$ be defined as before. Then
the set $S$ is called a Gr\"{o}bner-Shirshov basis  in  $RS\langle
X\rangle$ if any composition in $S$  is trivial modulo $S$.
\end{definition}

\begin{lemma}\label{15}
Let $S\subset RS\langle X\rangle$ be a set of monic polynomials and
$(u)_s$ an $S$-word. If any right multiplication composition in $S$
is trivial module $S$, then $(u)_s$ has a representation:
$$
(u)_s=\sum\limits_i \alpha_i [u_i]_{s_i}
$$
where each $\alpha_i\in k$, $[u_i]_{s_i}$ is a normal $S$-word and
$\overline{[u_i]_{s_i}}\leq \overline{(u)_s}$.
\end{lemma}

\noindent\textbf{Proof.} We use induction on $\overline{(u)_s}$. If $\overline{(u)_s}=\bar s$,  then $%
(u)_s=s$ and the result holds. Assume that $\overline{(u)_s}>\bar s$. Then $(u)_s=(v)_s[w]$ or $%
(u)_s=[w](v)_s$. We consider only the case $(u)_s=(v)_s[w]$. The
other case can be similarly proved.

By induction, we may assume that $(u)_s=[v]_s[w]$. If $[v]_s=s$,
then the result holds clearly because each right multiplication
composition in $S$ is trivial modulo $S$. Suppose that
$[v]_s=[v_1]_s[v_2]$ or $[v]_s=[v_1][v_2]_s$. We consider only the
case $[v]_s=[v_1]_s[v_2]$. The other case can be similarly proved.
If $[v_2]\leq [w]$, then $(u)_s=[v]_s[w]$ is a normal $S$-word and
we get the result. If $[v_2]> [w]$, then
$$
(u)_s=([v_1]_s[v_2])[w]=([v_1]_s[w])[v_2]+[v_1]_s([v_2][w])-[v_1]_s([w][v_2]).
$$
By induction, $[v_1]_s[w]=\sum\limits_j\beta_j[v_j]_{s_j}$, where
$\beta_j\in k$, $[v_j]_{s_j}$ is a normal $S$-word, and
$\overline{[v_j]_{s_j}}\leq \overline{[v_1]_s[w]}$. If
$\overline{[v_j]_{s_j}}= \overline{[v_1]_s[w]}$, then
$[v_j]_{s_j}[v_2]$ is a normal $S$-word since
$[v_j]_{\overline{s_j}}[v_2]=\overline{[v_j]_{s_j}}[v_2]=\overline{([v_1]_s[w])}[v_2]=
\overline{([v_1]_s[v_2])[w]}=\overline{(u)_s}$ by Corollary
\ref{ll5}. If $\overline{[v_j]_{s_j}}<\overline{[v_1]_s[w]}$, then
$\overline{[v_j]_{s_j}[v_2]}<\overline{([v_1]_s[w])[v_2]}=\overline{([v_1]_s[w])}[v_2]=\overline{(u)_s}$.
By Corollary \ref{ll5} again,
$$
 \overline{[v_1]_s([v_2][w])},
\overline{[v_1]_s([w][v_2])}<\overline{([v_1]_s[w])[v_2]}=\overline{(u)_s}.
$$
Now the result follows from the induction.\ \ $\square$

\begin{lemma}\label{16}
Let $[a_1s_1b_1],~[a_2s_2b_2]$ be normal $S$-words. If $S$ is a
Gr\"{o}bner-Shirshov basis in $RS\langle X\rangle$ and
$[w]=[a_1\overline{s_1}b_1]=[a_2\overline{s_2}b_2]$, then
\begin{equation*}
[a_1s_1b_1]\equiv [a_2s_2b_2] \ \ mod (S,[w]).
\end{equation*}
\end{lemma}

\noindent\textbf{Proof.} We have $w=a_1\bar s_1 b_1=a_2\bar s_2 b_2$
as associative words. There are two cases to consider.

Case 1. Suppose that  $\bar s_1$ and $\bar s_2$  are disjoint, say,
$|a_2|\geq |a_1|+|\bar s_1|$. Then, we can assume that
$$
a_2=a_1\bar s_1 c \ \ and  \ b_1=c\bar s_2 b_2
$$
for some $c\in X^*$. Thus, $ [w]=[a_1\bar s_1 c \bar s_2 b_2]. $
Now,
\begin{eqnarray*}
[a_1 s_1 b_1]-[a_2 s_2 b_2]&=&[a_1 s_1 c \bar s_2
b_2]-[a_1\bar s_1 c  s_2 b_2]\\
&=&[a_1 s_1 c \bar s_2 b_2]-(a_1 s_1 c  s_2 b_2)+(a_1 s_1 c  s_2
b_2)-[a_1\bar s_1 c  s_2 b_2]\\
&=&(a_1 s_1 c (\bar s_2 - s_2) b_2)+(a_1(s_1-\bar s_1) c  s_2 b_2).
\end{eqnarray*}
Since $[\overline{\overline{s_2}-s_2}]<\bar s_2$ and
$[\overline{s_1-\overline{s_1}}]<\bar s_1$, and by  Lemmas \ref{l4},
\ref{15}, we conclude that
$$
[a_1 s_1 b_1]-[a_2 s_2 b_2]=\sum\limits_i \alpha_i[u_is_iv_i]
$$
for some $\alpha_i\in k$, normal S-words $[u_is_iv_i]$ such that $
[u_i\bar s_1v_i]<[w]. $ Hence,
$$
[a_1s_1b_1]\equiv [a_2s_2b_2]\ mod (S,[w]).
$$

Case 2. Suppose that  $\bar s_1$  contains $\bar s_2$ as a subword.
We assume that
$$
\bar s_1=[a\bar s_2b], \ a_2=a_1a  \mbox{ and } b_2=bb_1, \mbox{
that is, } [w]=[a_1[a\bar s_2b]b_1]
$$
for the normal $S$-word $[a s_2 b]$. We have
\begin{eqnarray*}
[a_1 s_1 b_1]-[a_2 s_2 b_2]&=&[a_1 s_1 b_1]-[a_1 [a s_2 b] b_1]\\
&=&(a_1(s_1-[as_2b])b_1)\\
&=&(a_1(s_1,s_2)_{\bar s_1}b_1).
\end{eqnarray*}
 Since $S$ is a
Gr\"{o}bner-Shirshov basis, $(s_1,s_2)_{\bar
s_1}=\sum\limits_i\alpha_i[c_i s_i d_i]$ for some $\alpha_i\in k$,
normal S-words $[c_is_id_i]$ with each $[c_i\bar s_i d_i]<\bar s_1$.
By Lemmas \ref{l4}, \ref{15}, we have
\begin{eqnarray*}
&&[a_1 s_1 b_1]-[a_2 s_2 b_2]=(a_1(s_1,s_2)_{\bar s_1}b_1)\\
&=&\sum\limits_i\alpha_i(a_1[c_is_id_i]b_1)=\sum\limits_j\beta_j[a_js_jb_j]
\end{eqnarray*}
for some $\beta_j\in k $, normal S-words $[a_js_jb_j]$ with each
$[a_j\bar s_j b_j]<[w]=[a_1\bar {s_1}b_1]$. \\

Consequently, $ [a_1s_1b_1]\equiv [a_2s_2b_2]\ \  mod (S,[w]).  \ \
\ \ \square $

\begin{lemma}\label{l7}
Let $S\subset RS\langle X\rangle$ be a set of monic polynomials and
\ $Irr(S)=\{[u]\in W |[u]\ne [a\bar s b]\ a,b\in X^*,\ s\in S \mbox{
and } [as b] \mbox{ is a normal } S\mbox{-word}\}$. Then for any
$f\in RS\langle X\rangle$,
\begin{equation*}
f=\sum\limits_{[u_i]\leq \bar f }\alpha_i[u_i]+
\sum\limits_{[a_j\bar s_jb_j]\leq\bar f}\beta_j[a_js_jb_j],
\end{equation*}
where each $\alpha_i,\beta_j\in k, \ [u_i]\in Irr(S)$ and
$[a_js_jb_j]$ is a normal $S$-word.
\end{lemma}
{\bf Proof.} Let  $f=\sum\limits_{i}\alpha_{i}[u_{i}]\in{RS\langle
X\rangle}$ where $0\neq{\alpha_{i}\in{k}}$ and
$[u_{1}]>[u_{2}]>\cdots$. If $[u_1]\in{Irr(S)}$, then let
$f_{1}=f-\alpha_{1}[u_1]$. If $[u_1]\not\in{Irr(S)}$, then there
exist some $s\in{S}$ and $a_1,b_1\in{X^*}$ such that $\bar
f=[a_1\bar{s_1}b_1]$. Let $f_1=f-\alpha_1[a_1s_1b_1]$. In both
cases, we have $\bar{f_1}<\bar{f}$. Then the result follows by using
induction on $\bar{f}$. \ \ \ \ $\square$

\begin{theorem}\label{t1}
Let $S\subset RS\langle X\rangle$ be a nonempty set of monic
polynomials and the order $<$ be defined as before. Let $Id(S)$ be
the ideal of $RS\langle X\rangle$ generated by $S$. Then the
following statements are equivalent:
\begin{enumerate}
\item [(i)] $S$ is a Gr\"{o}bner-Shirshov basis in $RS\langle X\rangle$.

\item [(ii)] $f\in Id(S)\Rightarrow \bar f =[a\bar s b]$ for some $s\in S\
and\ a,b\in X^*$, where $[as b]$ is  a normal $S$-word.

\item [$(ii)'$] $f\in Id(S)\Rightarrow  f = \alpha_1[a_1s_1
b_1]+\alpha_2[a_2s_2b_2]+\cdots$, where $\alpha_i\in k, \ [a_1\bar
s_1b_1]>[a_2\bar s_2b_2]>\cdots$ and each $[a_is_i b_i]$ is  a
normal $S$-word.

\item [(iii)] $Irr(S)=\{[u]\in W |[u]\ne [a\bar s b]\ a,b\in X^*,\ s\in S \mbox{ and }
[as b] \mbox{ is a normal } S\mbox{-word}\}$ is a linear basis of
the algebra $RS\langle X | S\rangle=RS\langle X\rangle/Id(S)$.
\end{enumerate}
\end{theorem}

{\bf Proof.} $(i)\Rightarrow (ii)$. \ Let $S$ be a
Gr\"{o}bner-Shirshov basis and $0\neq f\in Id(S)$. We can also
assume, by Lemma \ref{15}, that
$$
f=\sum_{i=1}^n\alpha_i[a_is_ib_i]
$$
where each $\alpha_i\in k, \ a_i,b_i\in {X^*}, \ s_i\in S$ and $
[a_is_ib_i]$  a normal $S$-word. Let
$$
[w_i]=[a_i\bar s_ib_i], \
[w_1]=[w_2]=\cdots=[w_l]>[w_{l+1}]\geq\cdots
$$
We will use the induction on $l$ and $[w_1]$ to prove that
$\overline{f}=[a\overline{s}b]$ for some $s\in S \ \mbox{and} \
a,b\in {X^*}$.

If $l=1$, then $\overline{f}=\overline{[a_1s_1b_1]}=[a_1\bar
s_1b_1]$ and hence the result holds. Assume that $l\geq 2$. Then, by
Lemma \ref{16}, we have
$$
[a_1s_1b_1]\equiv[a_2s_2b_2] \ \ mod(S,[w_1]).
$$
Thus, if $\alpha_1+\alpha_2\neq 0$ or $l>2$, then the result holds.
For the case that $\alpha_1+\alpha_2= 0$ and $l=2$, we use the
induction
on $[w_1]$. Hence, the result follows.\\

$(ii)\Rightarrow (ii)'$. \ Assume that (ii) and $0\neq f\in Id(S)$.
Let $f=\alpha_1\overline{f}+\cdots$. Then, by (ii),
$\overline{f}=[a_1\bar s_1b_1]$. Therefore,
$$
f_1=f-\alpha_1[a_1s_1b_1], \ \overline{f_1}<\overline{f}, \ f_1\in
Id(S).
$$
Now, by using induction on $\overline{f}$, we have $(ii)'$.\\

$(ii)'\Rightarrow (ii)$. This part is clear.\\

$(ii)\Rightarrow(iii)$. Suppose that
$\sum\limits_{i}\alpha_i[u_i]=0$ in $RS\langle X | S\rangle$, where
$\alpha_i\in k$, $[u_i]\in {Irr(S)}$. It means that
$\sum\limits_{i}\alpha_i[u_i]\in{Id(S)}$ in ${RS\langle X\rangle}$.
Then all $\alpha_i$ must be equal to zero. Otherwise, we have
$\overline{\sum\limits_{i}\alpha_i[u_i]}=[u_j]\in{Irr(S)}$ for some
$j$ which contradicts (ii).

Now, for any $f\in{RS\langle X\rangle}$, by Lemma \ref{l7}, we have
$$
f=\sum\limits_{[u_i]\in Irr(S),\ [u_i]\leq \bar f }\alpha_i[u_i]+
\sum\limits_{[a_j\bar s_jb_j]\leq\bar f}\beta_j[a_js_jb_j].
$$
Hence, (iii) follows.\\

$(iii)\Rightarrow(i)$. Applying Lemma \ref{l7} to a composition of
elements of $S$, we get by $(iii)$ that any composition is trivial
because any composition belongs to $Id(S)$. So $S$ is a
Gr\"{o}bner-Shirshov basis. \ \ \ \ $\square$

\section{Gr\"{o}bner-Shirshov basis for universal enveloping
right-symmetric algebra of a Lie algebra}

The universal enveloping right-symmetric algebra of a Lie algebra is
defined in the paper \cite{Se94}. In this section, we give a
Gr\"{o}bner-Shirshov basis for such an algebra.

\begin{theorem}
Let $(\mathscr{L},[,])$ be a Lie algebra with a well ordered basis
$\{e_i|\ i\in I\}$. Let
$$
[e_i,e_j]=\sum\limits_{m}\alpha_{ij}^me_m
$$
where $\alpha_{ij}^m\in k$. We denote
$\sum\limits_{m}\alpha_{ij}^me_m$  by $\{e_ie_j\}$. Let
$$
U(\mathscr{L})=RS\langle \{e_i\}_I| \ e_ie_j-e_je_i=\{e_ie_j\}, \
i,j \in I\rangle
$$
be the universal enveloping right-symmetric algebra of
$\mathscr{L}$. Let
\begin{eqnarray*}
S&=&\{f_{ij}=e_ie_j-e_je_i-\{e_ie_j\},\ i,j \in I \ \mbox{ and } \
i>j \}.
\end{eqnarray*}
Then \begin{enumerate}
\item[(i)] \ $S$ is a Gr\"{o}bner-Shirshov basis in $RS\langle X\rangle$ where $X=\{e_i\}_I$.
\item[(ii)] (\cite{Se94}) \ $\mathscr{L}$ can be embedded into the ``universal enveloping
right-symmetric algebra" $U(\mathscr{L})$ as a vector space.
\end{enumerate}
\end{theorem}

\noindent\textbf{Proof.} (i). It is clear that $
\overline{f_{ij}}=e_ie_j \ (i>j)$. So, there exists a unique kind of
composition $f_{ij}e_k \ (i>j>k)$. Then, in $RS\langle X\rangle$, we
have
\begin{align*}
&f_{ij}e_k -f_{ik}e_j+f_{jk}e_i-e_if_{jk}+e_jf_{ik}-e_kf_{ij}-\sum_{m}
\alpha_{jk}^mf_{im}-\sum_{m}\alpha_{ij}^mf_{km}-\sum_{m}\alpha_{ik}^mf_{mj}\\
=&(e_ie_j)e_k-(e_je_i)e_k- \{e_ie_j\}e_k-(e_ie_k)e_j+(e_ke_i)e_j+\{e_ie_k\}e_j \\
&+f_{jk}e_i-e_if_{jk}+e_jf_{ik}-e_kf_{ij}-\sum_{m}\alpha_{jk}^mf_{im}-\sum_{m}
\alpha_{ij}^mf_{km}-\sum_{m}\alpha_{ik}^mf_{mj}\\
=&(e_ie_k)e_j+e_i(e_je_k)-e_i(e_ke_j)-(e_je_k)e_i-e_j(e_ie_k)+e_j(e_ke_i)-
\{e_ie_j\}e_k-(e_ie_k)e_j+(e_ke_j)e_i\\
&+e_k(e_ie_j)-e_k(e_je_i)+\{e_ie_k\}e_j+f_{jk}e_i-e_if_{jk}+e_jf_{ik}-e_kf_{ij}\\
& -
\sum_{m}\alpha_{jk}^mf_{im}-\sum_{m}\alpha_{ij}^mf_{km}-\sum_{m}\alpha_{ik}^mf_{mj}\\
=&-(e_je_k)e_i+(e_ke_j)e_i+e_i(e_je_k)-e_i(e_ke_j)+e_j(e_ke_i)-e_j(e_ie_k)-e_k(e_je_i)+e_k(e_ie_j)\\
&-
\{e_ie_j\}e_k+\{e_ie_k\}e_j+f_{jk}e_i-e_if_{jk}+e_jf_{ik}-e_kf_{ij}-\sum_{m}
\alpha_{jk}^mf_{im}-\sum_{m}\alpha_{ij}^mf_{km}-\sum_{m}\alpha_{ik}^mf_{mj}\\
=&-\{e_je_k\}e_i+e_i\{e_je_k\}-e_j\{e_ie_k\}+e_k\{e_ie_j\}-
\{e_ie_j\}e_k+\{e_ie_k\}e_j\\
&-\sum_{m}\alpha_{jk}^mf_{im}-\sum_{m}\alpha_{ij}^mf_{km}-\sum_{m}\alpha_{ik}^mf_{mj}\\
=&-\{\{e_ke_i\}e_j\}-\{\{e_ie_j\}e_k\}-\{\{e_je_k\}e_i\}\\
=&0\ \ (\mbox{by Jacobi identity}).
\end{align*}

By invoking $f_{mn}=-f_{nm}$, we have $f_{ij}e_k\equiv 0\ \ mod(S)$.
Therefore, $S$ is a Gr\"{o}bner-Shirshov basis for $U(\mathscr{L})$.

(ii). It follows from Theorem \ref{t1}. \ \ \ \ $\square$

\end{document}